\documentclass[10pt,twoside]{article}
\usepackage [latin1]{inputenc}
\usepackage[centertags]{amsmath}
\usepackage{amsfonts}
\usepackage{amssymb}
\usepackage{amsthm}
\usepackage{epsfig}
\usepackage{color}
\allowdisplaybreaks[4]

\newcommand{\abs}[1]{\lvert #1 \rvert}
\definecolor{c20}{rgb}{0.,0.7,0.}
\definecolor{c30}{rgb}{0.,0.,1.}
\definecolor{c40}{rgb}{1,0.1,0.7}
\definecolor{c50}{rgb}{1,0,0}
\definecolor{c60}{rgb}{1,0.9,0.1}

\def\zT#1{\textcolor{c50}{#1}}

\def\zT#1{#1}

\def\wHH{\mathcal{H}^2_{2H }}

\newcommand{\COM}[1]{}

\newcommand{\ABs}[1]{ \biggl \lvert #1 \biggr \rvert}
\newcommand{\R}{\mathbb{R}}

\newcommand{\inr}{\in \R}
\newcommand{\BQN}{\begin{eqnarray}}
\newcommand{\EQN}{\end{eqnarray}}
\newcommand{\BQNY}{\begin{eqnarray*}}
\newcommand{\EQNY}{\end{eqnarray*}}
\newtheorem{theo}{Theorem}[section]

\newtheorem{de}[theo]{Definition}
\newtheorem{lem}[theo]{Lemma}

\newtheorem{pro}[theo]{Proposition}

\newcommand{\BL}{\begin{lem}}
\newcommand{\EL}{\end{lem}}
\newcommand{\nelem}[1]{{Lemma \ref{#1}}}
\newcommand{\netheo}[1]{{Theorem \ref{#1}}}
\newcommand{\BT}{\begin{theo}}
\newcommand{\ET}{\end{theo}}
\newcommand{\BD}{\begin{de}}
\newcommand{\ED}{\end{de}}
\newcommand{\BP}{\begin{pro}}

\begin{document}

\newcommand{\EP}{\end{pro}}
\def\IF{\infty}
\newcommand{\E}[1]{\mathbb{E}\left \{ #1 \right \}}
\def\ZHT{Z(\tau,s)}
\def\ZHTY{\overline{Z}(\tau,s)}
\def\YHTY{\overline{\eta}(\tau,s)}
\def\MBT{M_{1/2}(T)}
\def\MHT{M_{H}(T)}
\def\MHTY{M^*_{H}(T)}
\def\MHTD{M_{H}^{(\delta)}(T)}
\def\MBTU{M_{H}(T_u)}
\def\MHTU{M_{H}(T_u)}

\title{\bf \Large Extremes of a type of locally stationary Gaussian random fields with applications to Shepp statistics\thanks{This work was supported by  Natural Science Foundation of Zhejiang Province of
China (No. LY18A010020) and National Science Foundation of China
(No. 11501250)}}
\author{{\small  Zhongquan Tan\footnote{E-mail address:  tzq728@163.com }\ \ \ Shengchao Zheng}\\
\\
{\small\it  College of Mathematics Physics and Information Engineering, Jiaxing University, Jiaxing 314001, PR China}\\
 }

\bigskip

\date{\today}
 \maketitle

{\bf Abstract:} Let $\{Z(\tau,s), (\tau,s)\in [a,b]\times[0,T]\}$ with some positive constants $a,b,T$ be
a centered Gaussian random field with variance function $\sigma^{2}(\tau,s)$ satisfying $\sigma^{2}(\tau,s)=\sigma^{2}(\tau)$.
 We firstly derive the exact tail asymptotics for
 the maximum
$M_{H}(T)=\max_{(\tau,s)\in[a,b]\times[0,T]}Z(\tau,s)/\sigma(\tau)$ up crossing some level $u$
with any fixed $0<a<b<\infty$ and $T>0$; and we further derive the extreme limit law for $M_{H}(T)$.
As applications of the main results, we derive the exact tail asymptotics and the extreme limit law
for Shepp statistics with stationary Gaussian process, fractional Brownian motion
and Gaussian integrated process as input.



\section{Introduction}
\def\TT{\mathcal{T}}
\def\IF{\infty}
\def\CCC{\mathbb{C}_*}
\def\XX{\{X(t),t\in [0,\IF)\}}

Let $B(t)$ be a standard Brownian motion and define the Shepp statistics as
$$\widehat{Z}(s) =\sup_{\tau\in[0,1]}B(s+\tau)-B(s),$$
for $s\geq 0$.
Since in various theoretical and applied problems, the Shepp statistics appears as the limit process
due to the central limit theorem, vast interest has been paid to the analysis of Shepp statistics,
see, e.g., Shepp (1971), Shepp and Slepian (1976), Cressie (1980), R\'{e}v\'{e}sz (1982), Deheuvels and Devroye (1987),
 Kabluchko and Munk (2008), Zholud (2008, 2009), Kabluchko (2011), Shklyaev (2011) and Tan (2015). A natural extension of the Shepp statistics is the following process,
$$\widehat{Z}_{H}(s) =\sup_{\tau\in[0,1]}B_{H}(s+\tau)-B_{H}(s),$$
where $B_H(t)$ is a standard fractional Brownian motion (fBm) with Hurst parameter $H\in (0,1)$.
This type of Shepp statistics have been extensively studied in
Hashorva and Tan (2013), D\c{e}bicki et al. (2015), Tan and Chen (2016) and Tan (2017).
Noting that the variance of the increment $B_{H}(s+\tau)-B_{H}(s)$ is $\tau^{2H}$, another type of natural extension of the Shepp statistics is
$$Z_{H}(s) =\sup_{\tau\in[a,b]}\frac{B_{H}(s+\tau)-B_{H}(s)}{\tau^H},\ \ 0<a<b<\infty.$$
Tan and Yang (2015) investigated the extreme $ M_{H}(T)= \sup_{0 \le s \le T}Z_{H}(s)$ and derived the following results.

\BT For $H\in(0, 1)$, $0<a<b<\infty$ and any $T>0$\zT{, it holds
that}
\begin{eqnarray}
\label{th2.1.1}
&&P\left(\max_{s\in[0,T]}Z_{H}(s) >u\right)=T\wHH \left(\frac{1}{2}\right)^{1/H}(1/a-1/b)u^{\frac{2}{H}}\Psi(u)(1+o(1)),
\end{eqnarray}
as $u\to \IF,$ and
\begin{eqnarray}
\label{f4}
\lim_{T\to \IF}  \max_{x\inr} \ABs{
P\left(a_{T}\left( \max_{s\in[0,T]}Z_{H}(s) -b_{T}\right)\leq x\right)-\exp(-e^{-x})}=0,
\end{eqnarray}
where
$$a_{T}=\sqrt{2 \ln T},\ \ b_{T}=a_T+a_T^{-1}\Bigl[ (\frac{1}{H}-\frac{1}{2}) \ln  \ln T+ \ln (2^{-\frac{1}{H}}
 \wHH (1/a-1/b)(2\pi)^{-1/2})\Bigl],$$
 $\Psi(u)$ and $\mathcal{H}_{2H}$ denote the tail distribution function of standard normal variable and Pickands constant (see definition in next section), respectively.
\ET

The Shepp statistics $Z_{H}(s)$ is a non-Gaussian random process, but it's extreme $M_{H}(T)$ can be transformed to
the extreme of a Gaussian random field $Z_{H}(\tau,s)$, i.e.,
$$M_{H}(T)=\max_{(\tau,s)\in[a,b]\times[0,T]}Z_{H}(\tau,s):=\max_{(\tau,s)\in[a,b]\times[0,T]}\frac{B_{H}(s+\tau)-B_{H}(s)}{\tau^H}.$$
It is easy to see from the Appendix that the random field $Z_{H}(\tau,s)$ is a locally stationary Gaussian random field.
It is natural to investigate the question
whether the above results still hold when $Z_{H}(\tau,s)$ is replaced by some more general locally stationary Gaussian random fields.
In this paper, we continue to study the limit properties of extremes of Shepp statistic and
pay our attention directly on the locally stationary Gaussian random field.
We first consider the tail asymptotics and the limit law of extremes for a type of locally stationary Gaussian random field.
Then we use the obtained results to derive the tail asymptotics and the limit law for the extremes of Shepp statistics with
stationary Gaussian process, fractional Brownian motion
and Gaussian integrated process as input, which extends Theorem 1.1.

The paper is organized as follows. Section 2 displays the main
results and Section 3 presents the applications.  The proofs of all results are \zT{postponed} to Sections 4.
In the Appendix, we present some useful definitions and
theorems for locally stationary Gaussian fields. In the following part of this paper,
let  $\Phi(\cdot)$  denote the probability distribution function of a standard normal variable and $\Psi(\cdot)=1-\Phi(\cdot)$.

\section{Main Results}

Suppose that $\{Z(\tau,s), (\tau,s)\in [a,b]\times[0,T]\}$ with some positive constants $a,b,T$ is a centered Gaussian random field with variance function and correlation function $\sigma^{2}(\tau,s)$ and
$r(\tau,s;\tau',s')$, respectively. Suppose the following assumptions hold.

{\bf Assumption A1:} there exists some positive  function $\sigma(\tau)$ which satisfies
$$\sigma(\tau,s)=\sigma(\tau),\ \  \forall (\tau,s)\in [a,b]\times[0,T].$$

{\bf Assumption A2:} there exist constant $\alpha\in (0,2]$ and continuous function $g(\tau)>0$ such that
$$r(\tau,s;\tau+\vartriangle_{\tau},s+\vartriangle_{ s})=1-g(\tau)(|\vartriangle_{ s}|^{\alpha}+|\vartriangle_{\tau}+\vartriangle_{ s}|^{\alpha})(1+o(1))$$
holds as $\vartriangle_{ s}\rightarrow0$ and $\vartriangle_{\tau}\rightarrow0$ and
further,
$$r(\tau,s;\tau',s')<1$$
for $(\tau,s)\neq(\tau',s')\in [a,b]\times[0,T]$.

{\bf Assumption A3:}  Assume that the function
$$\delta(v) := \sup\{|r(\tau,s;\tau',s')|, |s-s'|\ge v,s,s'\in [0,\infty), \tau,\tau'\in[a,b] \}$$ is such that
\BQN \label{BermD}
\lim_{v\to \IF} \delta(v)\ln v=r\in[0,\infty).
\EQN

The Pickands constant plays a crucial role in extreme value theory for Gaussian processes, which is defined by
$$\mathcal{H}_{\alpha}=\lim_{\lambda\rightarrow\infty} \lambda^{-1} \E{ \exp\left(\max_{t\in[0,\lambda]}
\sqrt{2}B_{\alpha/2}(t)-t^{\alpha}\right)} \in (0,\IF),$$ \zT{see e.g. Pickands
(1969), Piterbarg (1996).} Since we consider a random field, \zT{the
constant $\mathcal{H}_{\alpha}^2$ appears in our main results.} It is
well-known that $\mathcal{H}_{1}=1$ and $\mathcal{H}_{2}= 1/ \sqrt{
\pi}$.

\zT{Now we state our main results.}

\BT\label{Th:main21} Let $\{Z(\tau,s), (\tau,s)\in [a,b]\times[0,T]\}$ with some positive constants $a,b,T$ be a centered Gaussian random field with a.s. continuous sample
paths. Suppose that assumptions {\bf A1-A2} are satisfied with the parameters mentioned therein.
Then
\begin{eqnarray}
\label{th2.1.1}
&&P\left(\max_{(\tau,s)\in[a,b]\times[0,T]} Z(\tau,s)/\sigma(\tau)  >u\right)=T\mathcal{H}_{\alpha}^{2} \int_{a}^{b}(g(t))^{\frac{2}{\alpha}}dtu^{\frac{4}{\alpha}}\Psi(u)(1+o(1)),
\end{eqnarray}
as $u\to \IF.$ \ET

In the following part of the paper, let $\mathcal{N}$ be a standard normal random variable.

\BT\label{Th:main22} Let $\{Z(\tau,s), (\tau,s)\in [a,b]\times[0,T]\}$ with some positive constants $a,b$ be a centered Gaussian random field with a.s. continuous sample paths. Suppose that assumptions {\bf A1-A3} are satisfied with the parameters mentioned therein.
In addition, assume that $\{Z(\tau,s), (\tau,s)\in [a,b]\times[0,T]\}$ is homogeneous with respect to the second factor $s$.
Then
\begin{eqnarray}\label{f4}
\lim_{T\to \IF}  \max_{x\inr} \ABs{
P\left(a_{T}\left( \max_{(\tau,s)\in[a,b]\times[0,T]}Z(\tau,s)/\sigma(\tau)  -b_{T}\right)\leq x\right)-E\exp(-e^{-x-r+\sqrt{2r}\mathcal{N}})}=0,
\end{eqnarray}
 and
$$a_{T}=\sqrt{2 \ln T},\ \ b_{T}=a_T+a_T^{-1}\Bigl[ (\frac{2}{\alpha}-\frac{1}{2}) \ln  \ln T+ \ln (
 \mathcal{H}_{\alpha}^{2} \int_{a}^{b}(g(t))^{\frac{2}{\alpha}}dt(2\pi)^{-1/2})\Bigl].$$\ET

In the literature assumption {\bf A3} with $r=0$ is referred to as the weak
dependence or the Berman's condition and consequently the random field $Z$ is called a weakly dependent. In
analogy, the random field $Z$ with correlation function satisfying assumption {\bf A3} with $r>0$  is called strongly dependent.

{\bf Remark 2.1:} i). Theorem 2.1 and 2.2 extend the main results of Tan and Yang (2015) from fBm field to general locally stationary Gaussian random fields.
Specially, we derived the strongly dependent case. \\
ii). It is worth mentioning the work of D\c{e}bicki et al. (2015) and Tan (2017), where the tail asymptotics and the extreme limit theorem
for a type of non-homogeneous Gaussian random field are derived, respectively. The method of proofs used in this paper are different to that of the aforementioned papers.

We end this section with an example which satisfies all of the cases of assumptions {\bf A1-A3}.

{\bf Example 2.1:}  Consider a Gaussian random field defined as
$$
Z(\tau,s)=\frac{1}{\sqrt{2}}(Y(\tau+s)+X(s))\sigma(\tau),\quad (\tau,s)\in [a,b]\times[0,T],
$$
where  $X,Y$ are two independent centered stationary Gaussian processes with  the same covariance functions $r$  satisfying as $t\to 0$
$$
r(t)=1-a|t|^{\alpha}(1+o(1)),
$$
for some constants $a>0,\alpha\in(0,2]$. Here $\sigma(t)$ is a positive function. Further, assume that
\BQNY
 r(t)<1,\ \forall\ t\neq 0,\ \ \ \lim_{t\rightarrow\infty}r(t)\ln t=r\in[0,\infty)
\EQNY
 It follows that the assumptions {\bf A1-A3} are satisfied by $\{Z(\tau,s), (\tau,s)\in [a,b]\times[0,T]\}$.

%\BT\label{Th:main2}  \ET
%\textbf{Remark 2.4}.

\section{Applications}

Throughout this section, let $\{X(t),t\geq0\}$ be a centered Gaussian process and define general Shepp statistics as
$$Y(s)=\max_{\tau\in[a,b]}Y(\tau,s)=\max_{\tau\in[a,b]}\frac{X(s+\tau)-X(s)}{\sqrt{E(X(s+\tau)-X(s))^{2}}},$$
for some fixed $a,b,T>0$.
Applying Theorems \ref{Th:main21} and \ref{Th:main22}, we study the limit properties of extremes of Shepp statistics for a more general Gaussian process $X(t)$, which is a stationary Gaussian process or non-stationary Gaussian process with stationary increments.

First, let $\{X(t),t\geq0\}$ be a stationary Gaussian process with mean 0, variance 1.
Suppose the covariance function $r_{X}$ of $\{X(t),t\geq0\}$ satisfies the following conditions:

{\bf Assumption B1:} there exist positive constant $a_{1}$ and $\alpha\in(0,2)$ such that
$r_{X}(t)=1-a_{1}|t|^{\alpha}(1+o(1))$, as $t\rightarrow 0;$ $r_{X}(t)\in C([0,\infty))$ and $r_{X}(s)<1$ for $s>0$.

{\bf Assumption B2:} $r_{X}(t)$ is twice continuously differentiable on $[\lambda,\infty)$ for some $\lambda>0$ and $\lim_{t\rightarrow\infty}\ddot{r}_{X}(t)\ln t=r\in[0,+\infty)$.

\BP\label{Pro:main3:1} Let $Y(\tau,s)$ be defined as above. Suppose that $r_{X}(t)$ satisfies Assumption ${\bf B1}$. Then we have
\begin{eqnarray}
\label{P3.1.1}
&&P\left(\max_{(\tau,s)\in[a,b]\times[0,T]} Y(\tau,s) >u\right)=T\mathcal{H}_{\alpha}^{2} (a_{1}/2)^{2/\alpha}\int_{a}^{b}(1-r_{X}(t))^{-2/\alpha}dtu^{\frac{4}{\alpha}}\Psi(u)(1+o(1)),
\end{eqnarray}
as $u\to \IF.$
Furthermore, if Assumption ${\bf B2}$ holds, then
\BQN \label{thm3:1}
\lim_{T\to \infty}\sup_{x\in \mathbb{R}}
\ABs{P\left(a_{T}\left(\sup_{(\tau,s)\in[a,b]\times[0,T]}Y(\tau,s) -b_{T}\right)\leq x\right)-E\exp(-e^{-x-r+\sqrt{2r}\mathcal{N}})}=0,
\EQN
where  $$a_{T}=\sqrt{2 \ln T},\ \ b_{T}=a_T+a_T^{-1}\Bigl[ (\frac{2}{\alpha}-\frac{1}{2}) \ln  \ln T+ \ln ((a_{1}/2)^{\frac{2}{\alpha}}
 \mathcal{H}_{\alpha}^{2} \int_{a}^{b}(1-r_{X}(t))^{-2/\alpha}dt(2\pi)^{-1/2})\Bigl].$$
\EP

{\bf Example 3.1:} There are many types of stationary Gaussian processes such as the fractional Ornstein-Uhlenbeck process with covariance function $r_{X}(t)=e^{-|t|^{\alpha}}$ and the generalized Cauchy model with covariance function $r_{X}(t)=(1+|t|^{\alpha})^{-\beta}$ with $\alpha\in(0,2)$ and $\beta>0$ which satisfy the conditions of Proposition \ref{Pro:main3:1} with $r=0$.

Second, let $\{X(t),t\geq0\}$ be a centered non-stationary Gaussian process with stationary increment and variance function $\sigma_{X}^{2}(t)$, a.s. continuous sample paths. Recall that $X(t)$ is said to have stationary increments if the law of the process $\{X(t+t_{0})-X(t_{0}),t\in \mathbb{R}\}$ does not depend on the choice of $t_{0}$.
Suppose that the variance function $\sigma_{X}^{2}(t)$ of $\{X(t),t\geq0\}$ satisfies the following conditions:

{\bf Assumption C1:} $\sigma_{X}^{2}(t)$ is twice continuously differentiable on $[\lambda,\infty)$ for $\lambda>0$  and further
$\sigma_{X}^{2}(t)=a_{2}|t|^{\alpha}(1+o(1))$, as $t\rightarrow 0$
holds for some $\alpha\in(0,2], a_{2}>0$.

{\bf Assumption C2:} $\lim_{t\rightarrow\infty}\ddot{\sigma}_{X}^{2}(t)\ln t\rightarrow r\in[0,+\infty)$.

\BP\label{pro:main3:1} Let $Y(s,t)$ be defined as above. Suppose that $\sigma_{X}(t)$ satisfies Assumption ${\bf C1}$. We have
 $$P\left(\sup_{(\tau,s)\in[a,b]\times[0,T]}Y(\tau,s)>u\right)=T\mathcal{H}_{\alpha}^{2} (a_{2}/2)^{2/\alpha}\int_{a}^{b}(\sigma_{X}(t))^{-\frac{4}{\alpha}}dtu^{\frac{4}{\alpha}}\Psi(u)(1+o(1)),$$
 as $u\rightarrow\infty$.
Furthermore, if Assumption ${\bf C2}$ holds, then
\BQN \label{thm3:1}
\lim_{T\to \infty}\sup_{x\in \mathbb{R}}
\ABs{P\left(a_{T}\left(\sup_{(\tau,s)\in[a,b]\times[0,T]}Y(\tau,s) -b_{T}\right)\leq x\right)-E\exp(-e^{-x-r+\sqrt{2r}\mathcal{N}})}=0,
\EQN
where  $$a_{T}=\sqrt{2 \ln T},\ \ b_{T}=a_T+a_T^{-1}\Bigl[ (\frac{2}{\alpha}-\frac{1}{2}) \ln  \ln T+ \ln ((a_{2}/2)^{\frac{2}{\alpha}}
 \mathcal{H}_{\alpha}^{2} \int_{a}^{b}(\sigma_{X}(t))^{-\frac{4}{\alpha}}dt(2\pi)^{-1/2})\Bigl].$$
\EP

The following two examples satisfy conditions of Proposition \ref{pro:main3:1}. The first example extends the main results of Tan and Yang (2015).

{\bf Example 3.2:}  Let $B_{H_{i}}(t)$, $i=1,2,\ldots,n$ be a sequence of independent fBms with Hurst index $H_{i}\in(0,1)$ and $\lambda_{i}$ be a positive
sequence satisfying $\sum_{i=1}^{n}\lambda_{i}^{2}=1$.
Since $\lambda_{1}B_{H_{1}}(t)+\lambda_{2}B_{H_{2}}(t)=^{d}\sqrt{\lambda_{1}^{2}+\lambda_{2}^{2}}B_{H}(t)$ for $H=H_{1}=H_{2}$ , we suppose that
$$H:=H_{1}<H_{2}<\cdots<H_{n}.$$
Let
$X(t)=\sum_{i=1}^{n}\lambda_{i}B_{H_{i}}(t)$
and $Y(\tau,s)$ be defined as above. We have
$$P\left(\sup_{(\tau,s)\in[a,b]\times[0,T]}Y(\tau,s) >u\right)=T\mathcal{H}_{2H}^{2} (1/2)^{1/H}\int_{a}^{b}\frac{1}{(\sum_{i=1}^{n}\lambda_{i}^{2}|t|^{2H_{i}})^{1/H}}dtu^{\frac{2}{H}}\Psi(u)(1+o(1)),$$
as $u\rightarrow\infty$,
and
\BQN \label{thm3:2}
\lim_{T\to \infty}\sup_{x\in \mathbb{R}}
\ABs{P\left(a_{T}\left(\sup_{(\tau,s)\in[a,b]\times[0,T]}Y(\tau,s) -b_{T}\right)\leq x\right)-\exp\{-e^{-x}\}}=0,
\EQN
where  $$a_{T}=\sqrt{2 \ln T},\ \ b_{T}=a_T+a_T^{-1}\Bigl[ (\frac{1}{H}-\frac{1}{2}) \ln  \ln T+ \ln ((1/2)^{\frac{1}{H}}
 \mathcal{H}_{2H}^{2} \int_{a}^{b}\frac{1}{(\sum_{i=1}^{n}\lambda_{i}^{2}|t|^{2H_{i}})^{1/H}}dt(2\pi)^{-1/2})\Bigl].$$

The next example considers the Gaussian integrated process. For related studies, we refer to D\c{e}bicki (2002) and H\"{u}sler and Piterbarg (2004).

{\bf Example 3.3:} Let $\{\zeta(t),t\geq0\}$ be a centered stationary Gaussian process with variance one and
suppose the covariance function $r_{\zeta}(t)$ of $\{\zeta(t),t\geq0\}$ satisfying the following conditions:

{\bf Assumption D1:} $r_{\zeta}(t)\in C([0,\infty))$ and $\int_{0}^{t}r_{\zeta}(s)ds>0$ for $t\in(0,T]$; $r_{\zeta}(t)=1-t^{\theta}(1+o(1))$ as $t\rightarrow 0^{+}$ with $\theta\in(0,2]$;

{\bf Assumption D2:} $\lim_{t\rightarrow\infty}r_{\zeta}(t)\ln t=r\in[0,+\infty)$.

Define Gaussian integrated processes as
$X(t)=\int_{0}^{t}\zeta(s)ds$
and let $Y(\tau,s)$ be defined as before.
If Assumption {\bf D1} are satisfied, we have for some constant $T>0$
$$P\left(\sup_{(\tau,s)\in[a,b]\times[0,T]}Y(\tau,s) >u\right)=T\pi^{-1} (1/4)\int_{a}^{b}\left(\int_{0}^{t}(t-s)r_{\varsigma}(s)ds\right)^{-1}dtu^{2}\Psi(u)(1+o(1))$$
as $u\rightarrow\infty$.
If further Assumption {\bf D2} holds, we have
\BQN \label{thm3:2}
\lim_{T\to \infty}\sup_{x\in \mathbb{R}}
\ABs{P\left(a_{T}\left(\sup_{(\tau,s)\in[a,b]\times[0,T]}Z(\tau,s) -b_{T}\right)\leq x\right)-E\exp(-e^{-x-r+\sqrt{2r}\mathcal{N}})}=0,
\EQN
where  $$a_{T}=\sqrt{2 \ln T},\ \ b_{T}=a_T+a_T^{-1}\Bigl[ \frac{1}{2} \ln  \ln T+ \ln ((1/4\pi)
\int_{a}^{b}\left(\int_{0}^{t}(t-s)r_{\varsigma}(s)ds\right)^{-1}dt(2\pi)^{-1/2})\Bigl].$$

\section{Proofs}

In this section, we give \zT{the} detailed proofs of Theorems
2.1 and 2.2 and Propositions 3.1 and 3.2.

\subsection{Proof of Theorem 2.1}
As in the proof of Theorem 2.1 of Tan and Yang (2015), we use \netheo{Theorem
4.1} \zT{in} Appendix  to prove the theorem.

\textbf{Proof:} By the definition, it
follows that the standardised random field $
Z(\tau,s)/\sigma(\tau)$ is locally stationary (see Appendix 5 for
details). The local structure is given
by
$$C_{(\tau,s)}(\vartheta,v)=g(\tau)(|\vartheta+v|^{\alpha}+|v|^{\alpha}),$$
where $ (\tau,s)\in K=\{(\tau,s):\tau\in [a,b],s\in[0,T]\}$
and $(\vartheta,v)\in \mathbb{R}^{2}$. Thus, the tangent process of
$Z(\tau,s)/\sigma(\tau)$ is given by
$$
Y_{(\tau,s)}(\vartheta,v)=g(\tau)\chi(\vartheta,v),$$ where
$\chi(\vartheta,v)=B_{\alpha/2}(\vartheta+v)+\tilde{B}_{\alpha/2}(v) - |\vartheta+v|^{\alpha}- |v|^{\alpha}$ with $B_{\alpha/2}$
and $\tilde B_{\alpha/2}$ two independent standard fBm's. By \netheo{Theorem
4.1} \zT{in} Appendix
$$P\left(\max_{(\tau,s)\in[a,b]\times[0,T]} Z(\tau,s)/\sigma(\tau)  >u\right)
=\left(\int_{K}h(\tau,s)d\tau ds
\right)u^{\frac{4}{\alpha}}\Psi(u)(1+ o(1)), \quad u\to
\IF,$$ where $h(\tau,s)$ is defined by
$$
h(\tau,s)=\lim_{\mathcal{T}\rightarrow\infty}\frac{1}{\mathcal{T}^{2}}\mathbb{E}
\left\{\exp\left(
\max_{(\vartheta,v)\in[0,\mathcal{T}]\times[0,\mathcal{T}]}
g(\tau)\chi(\vartheta,v)\right)\right\}.$$ Since the
field $\chi(\vartheta,v)+ \abs{\vartheta+v}^{\alpha}+ \abs{v}^{\alpha}$ is $\alpha$ self-similar\zT{, we have} \BQNY
\lefteqn{\int_{\tau\in [a,b], s\in [0,T]} h(\tau,s)\,d\tau ds }\\
&=&\int_{\tau\in [a,b], s\in [0,T]}
\lim_{\mathcal{T}\rightarrow\infty}\frac{1}{\mathcal{T}^{2}}\mathbb{E}
\left\{\exp\left( \max_{(\vartheta,v)\in[0,\mathcal{T}]\times[0,\mathcal{T}]} g(\tau)\chi(\vartheta,v)\right)\right\}d\tau ds \\
&=&  \int_{\tau\in [a,b], s\in [0,T]}
(g(\tau))^{2/\alpha}
\lim_{\mathcal{T}\rightarrow\infty}\frac{1}{\mathcal{T}^{2}}\mathbb{E}
\left\{\exp\left( \max_{(\vartheta,v)\in[0,\mathcal{T}]\times[0,\mathcal{T}]} \chi(\vartheta,v)\right)\right\}d\tau ds\\
&=& \mathcal{H}_{\alpha}^{2} \int_{\tau\in [a,b], s\in [0,T]}(g(\tau))^{2/\alpha}d\tau ds\\
&=& \mathcal{H}_{\alpha}^{2}T\int_{a}^{b}(g(\tau))^{2/\alpha}d\tau. \EQNY
Consequently, %Thus, we have
\BQNY P\left(\max_{(\tau,s)\in[a,b]\times[0,T]} Z(\tau,s)/\sigma(\tau)  >u\right)
&=&\mathcal{H}_{\alpha}^{2}T\int_{a}^{b}(g(\tau))^{2/\alpha}d\tau
u^{\frac{4}{\alpha}}\Psi(u)(1+o(1)), \quad u\rightarrow\infty, \EQNY
hence the first claim follows.  \hfill$\Box$

\subsection{Proof of Theorem 2.2}
In the following, for simplicity, let $ \ZHTY= \ZHT/\sigma(\tau)$,
 $u=u_{T}(x)=a_{T}^{-1}x+b_{T}$ and write $$w(u)=\mathcal{H}_{\alpha}^{2}\int_{a}^{b}(g(t))^{\frac{2}{\alpha}}dtu^{\frac{4}{\alpha}}\Psi(u).$$
 It is easy to check that $Tw(u_{T})\rightarrow e^{-x}$ as $T\rightarrow\infty$, by the definitions of $a_{T}$ and $b_{T}$.
First, we divide the interval [0,T] into intervals with constant length $L$ alternating with shorter intervals of length $\delta<L$. Define
$$I_{k}=[a,b]\times [(k-1)L,kL-\delta), \quad I_{k}^{*}=[a,b]\times [kL-\delta,kL),$$
then we have with $\label{KT} K_{T}:=\left[\frac{T}{L}\right]\in
\mathbb{N},$
 \BQNY [a,b]\times [0,T]
=\bigcup_{k=1}^{K_{T}}(I_{k}\cup I_{k}^{*})\cup I_{K_{T}+1}, \text{
where }I_{K_{T}+1}=[a,b]\times [K_{T}L,T], \EQNY \zT{which
implies} the length of the last interval $|I_{K_{T}+1}|\leq L$.
Obviously, we can apply Theorem
2.1 for these short and long intervals, respectively, since $\delta$ and $L$ are independent of $u$.

Next, we will construct a Gaussian random field to approximate $Z(\tau,s)$.
Let $Z_{i}(\tau,s)$, $i=1,2,\ldots$ be independent copy of $Z(\tau,s)$ and let $\xi(\tau,s)$ be such that
$\xi(\tau,s)=Z_{j}(\tau,s)$ for $(\tau,s)\in I_{j}$. Let $\rho(T)=r/(\ln T)$, and define
$$\eta(\tau,s)=(1-\rho(T))^{1/2}\xi(\tau,s)+(\rho(T))^{1/2}\sigma(\tau)\mathcal{N},\ \ (\tau,s)\in \cup_{j=1}^{K_{T}}I_{j},$$
where $\mathcal{N}$ is a standard Gaussian random variable which is independent of   $\xi(\tau,s)$.
Denote by $\varrho(\tau,s;\tau',s')$ the correlation function of $\eta(\tau,s)$. It is easy to check that
\[
  \varrho(\tau,s;\tau',s')=\left\{
 \begin{array}{cc}
  {r(\tau,s;\tau',s')+(1-r(\tau,s;\tau',s'))\rho(T)},    & (\tau,s)\in I_{j}, (\tau',s')\in I_{i}, i=j,\\
  {\rho(T)},    & (\tau,s)\in I_{j}, (\tau',s')\in I_{i}, i\neq j.
 \end{array}
  \right.
\]
Note that $Var(\eta(\tau,s))=Var(Z(\tau,s))=\sigma^{2}(\tau)$. In the sequel, $C$ shall denote
positive constant whose values may vary from place to place.

\BL \label{Lemma 3.2} \zT{From} the definitions of $I_{k}$, $k\geq
1$, it follows that \zT{as} $T\rightarrow\infty$ and
$\delta\downarrow0$
\begin{eqnarray}
\label{eq4.3}
\left|P\left( \max_{\tau\in[a,b]\atop s\in[0,T]} \ZHTY \leq u\right)-
P\left( \max_{(\tau,s)\in \cup_{k=1}^{K_{T}} I_{k}} \ZHTY \leq u\right)\right|\rightarrow0.
\end{eqnarray} \EL

\textbf{Proof:}  By the definitions of $I_{k}$ and $I_{k}^{*}$, we
\zT{rewrite}
\begin{eqnarray*}
P\left( \max_{\tau\in[a,b]\atop s\in[0,T]} \ZHTY \leq u\right)
&=&P\left( \max_{(\tau, s)\in\cup_{k=1}^{K_{T}}(I_{k}\cup I_{k}^{*})\cup I_{K_{T}+1}} \ZHTY \leq u\right).
\end{eqnarray*}
Thus, in order to prove (\ref{eq4.3}), it suffices to show that
\begin{eqnarray}
\label{eq4.31}
T_{u}:= \left|P\left(\max_{(\tau, s)\in\cup_{k=1}^{K_{T}}(I_{k}\cup I_{k}^{*})\cup I_{K_{T}+1}} \ZHTY \leq u\right)-
P\left( \max_{(\tau,s)\in \cup_{k=1}^{K_{T}} I_{k}} \ZHTY \leq u\right)\right|\rightarrow 0,
\end{eqnarray}
 as $T\rightarrow\infty$.
\zT{Obviously, for sufficiently large} $u$,
\begin{eqnarray}
\label{eq4.32}
T_u \le P\left( \max_{(\tau, s)\in \cup_{k=1}^{K_{T}} I_{k}^{*}\cup I_{K_{T}+1}} \ZHTY > u\right)
&\leq&\sum_{k=1}^{K_{T}}P\left(\max_{(\tau, s)\in I_{k}^{*}} \ZHTY > u\right)\nonumber\\
&+&P\left(\max_{(\tau, s)\in I_{K_{T}+1}} \ZHTY > u\right).
\end{eqnarray}
\zT{By Theorem 2.1,} the right-hand side of (\ref{eq4.32}) is
bounded by
\begin{eqnarray}
\label{eq4.33}
C\delta K_{T}w(u)+CLw(u)
&=& C\delta\frac{K_{T}}{T}Tw(u)+C\frac{L}{T}Tw(u)\nonumber\\
&=&C\frac{\delta}{L}+C\frac{L}{T}
\end{eqnarray}
where in the last step we use the fact that $Tw(u)=O(1)$, as $T\rightarrow\infty$. Since $L$ is a positive
constant, we conclude that the right-hand side of (\ref{eq4.33})
tends to 0 as
$T\rightarrow\infty$ and $\delta\downarrow0$. Thus, (\ref{eq4.31})
follows, and \zT{this completes the proof of the lemma}. \hfill$\Box$

\COM{Next, we shall apply Berman's
inequality to prove that the maximum  on the intervals $I_{k}$ are
asymptotically independent. Since the Berman's inequality \zT{ holds
only for the} sequences of Gaussian random variables,}

Next, we will approximate the continuous time maximum of $\ZHTY$ by a discrete one, so introduce
the following grids points. For some small $d>0$ and any $u$,
we define a family of grid points as follows. Let
$$q=q(u)=du^{-\frac{2}{\alpha}}$$ and define the grid of points \BQN \label{tauj}
s_{k,l}=(k-1)L+lq\ \ \mbox{and}\ \ \tau_{j}=b-jq, \EQN with
$(\tau_{j}, s_{k,l})\in I_{k}$ for integers $j\in \mathbb{Z}^{+}$,
$l\geq 0$, $k\geq 1$. These grid points are \zT{simply denoted} by
$(\tau, s)\in I_{k}\cap \mathcal{R}$ for fixed $k$, without
mentioning the dependence on $u$.

\BL \label{Lemma 3.3} Let $\YHTY=\eta(\tau,s)/\sqrt{Var(\eta(\tau,s))}$. \zT{It holds that} as $T
\to \IF$ and $d\downarrow 0$
\begin{eqnarray}
\label{eq4.10} 0\leq P\left(\max_{(\tau,s)\in \cup I_{k}\cap
\mathcal{R}} \ZHTY \leq u\right) -P\left( \max_{(\tau,s)\in \cup
I_{k}} \ZHTY \leq u\right)\rightarrow0,
\end{eqnarray}
\zT{as well as}
\begin{eqnarray}
\label{eq4.11}
0\leq P\left(\max_{(\tau,s)\in \cup I_{k}\cap
\mathcal{R}} \YHTY \leq u\right) -P\left( \max_{(\tau,s)\in \cup
I_{k}} \YHTY \leq u\right)\rightarrow0.
\end{eqnarray}
\EL

\textbf{Proof:} Applying Theorem \ref{Theorem 4.2} in  Appendix, repeating the proof of Theorem 2.1, we have for any $k$
$$P\left( \max_{(\tau,s)\in I_{k}\cap \mathcal{R}}\ZHTY> u\right)
   =L\mathcal{H}_{\alpha,d}^{2} \int_{a}^{b}(g(t))^{\frac{2}{\alpha}}dtu^{\frac{4}{\alpha}}\Psi(u)(1+o(1))(1+o(1))$$
as $T\rightarrow\infty$, where $\mathcal{H}_{\alpha,d}$ is a type of Pickands constant defined by
$$\mathcal{H}_{\alpha,d}=\lim_{\lambda\rightarrow\infty} \lambda^{-1} \E{ \exp\left(\max_{kd\in[0,\lambda]}
\sqrt{2}B_{\alpha/2}(kd)-(kd)^{\alpha}\right)} \in (0,\IF),$$
satisfying $\lim_{d\downarrow0}\mathcal{H}_{\alpha,d}=\mathcal{H}_{\alpha}$.
By Theorem 2.1, we have for any $k$
$$P\left( \max_{(\tau,s)\in I_{k}}\ZHTY> u\right)
   =L\mathcal{H}_{\alpha}^{2}\int_{a}^{b}(g(t))^{\frac{2}{\alpha}}dtu^{\frac{4}{\alpha}}\Psi(u)(1+o(1))(1+o(1))$$
as $T\rightarrow\infty$. \zT{Thus, we have }
\begin{eqnarray*}
0&\leq& P\left(\max_{(\tau,s)\in \cup I_{k}\cap \mathcal{R}} \ZHTY \leq u\right)
    -P\left(\max_{(\tau,s)\in \cup I_{k}} \ZHTY \leq u\right)\\
&\leq& \sum_{k=1}^{K_{T}}\left(P\left(\max_{(\tau,s)\in  I_{k}\cap \mathcal{R}} \ZHTY \leq u\right)
    -P\left(\max_{(\tau,s)\in  I_{k}} \ZHTY \leq u\right)\right)\\
&\leq&LK_{T}(\mathcal{H}_{\alpha,d}^{2}-\mathcal{H}_{\alpha}^{2})\int_{a}^{b}(g(t))^{\frac{2}{\alpha}}dtu^{\frac{4}{\alpha}}\Psi(u)
=:P_{T}.
\end{eqnarray*}
By the definition of $u=u_{T}(x)=a_{T}^{-1}x+b_{T}$, it is easy to check that
$$T\mathcal{H}_{\alpha}^{2}\int_{a}^{b}(g(t))^{\frac{2}{\alpha}}dtu^{\frac{4}{\alpha}}\Psi(u)=e^{-x}(1+o(1)),$$
as $T\rightarrow\infty$. Thus, we have
$$P_{T}\leq C (\mathcal{H}_{\alpha,d}^{2}/\mathcal{H}_{\alpha}^{2}-1)\rightarrow 0,$$
as $d\rightarrow 0$, which \zT{leads to} (\ref{eq4.10}). The claim
(\ref{eq4.11}) can be proven according to the similar arguments.
\hfill$\Box$

\BL \label{Lemma 3.4} Let $\YHTY=\eta(\tau,s)/\sqrt{Var(\eta(\tau,s))}$
\zT{, it holds that} as $T\rightarrow\infty,$
\begin{eqnarray}
\label{eq4.12}
\left|P\left(\max_{(\tau,s)\in \cup I_{k}\cap \mathcal{R}} \ZHTY \leq u\right)
-P\left(\max_{(\tau,s)\in \cup I_{k}\cap \mathcal{R}} \YHTY \leq u\right)\right|\rightarrow0,
\end{eqnarray}
uniformly for $d>0$. \EL

\textbf{Proof:} Applying Berman's inequality (see e.g. Piterbarg (1996)), we have
\begin{eqnarray}\label{tt1}
&&\left|P\left(\max_{(\tau,s)\in \cup I_{k}\cap \mathcal{R}} \ZHTY \leq u\right)
-P\left(\max_{(\tau,s)\in \cup I_{k}\cap \mathcal{R}} \YHTY \leq u\right)\right|\nonumber\\
&&\leq C\sum_{(\tau_{j},s_{k,l}),(\tau_{j'},s_{k',l'})\in \cup I_{i}}
|r(\tau_{j},s_{k,l};\tau_{j'},s_{k',l'})-\varrho(\tau_{j},s_{k,l};\tau_{j'},s_{k',l'})|\nonumber\\
&&\ \ \ \ \ \times\int_{0}^{1}\frac{1}{\sqrt{1-r_{h}(\tau_{j},s_{k,l};\tau_{j'},s_{k',l'})}}\exp\left(-\frac{u^{2}}{1+r_{h}(\tau_{j},s_{k,l};\tau_{j'},s_{k',l'})}\right)dh\nonumber\\
&&=:\sum_{(\tau_{j},s_{k,l}),(\tau_{j'},s_{k',l'})\in I_{i},\atop 1\leq i\leq K_{T}}M(u)
+\sum_{(\tau_{j},s_{k,l})\in I_{i},(\tau_{j'},s_{k',l'})\in I_{i'},\atop 1\leq i\neq i'\leq K_{T}}M(u),
\end{eqnarray}
where
$$r_{h}(\tau_{j},s_{k,l};\tau_{j'},s_{k',l'})=hr(\tau_{j},s_{k,l};\tau_{j'},s_{k',l'})+(1-h)\varrho(\tau_{j},s_{k,l};\tau_{j'},s_{k',l'}),\ \ h\in[0,1].$$

Next, we estimate the upper bound for the first sum in (\ref{tt1}). Let
$$\varpi(\tau_{j},s_{k,l};\tau_{j'},s_{k',l'})=\max\{|r(\tau_{j},s_{k,l};\tau_{j'},s_{k',l'})|,|\varrho(\tau_{j},s_{k,l};\tau_{j'},s_{k',l'})|\}$$
and $$\vartheta(t)=\sup_{\{t\leq |s_{k,l}-s_{k',l'}|\leq T\}\cup\{t\leq |\tau_{j}-\tau_{j'}|\leq (b-a)\}}\varpi(\tau_{j},s_{k,l};\tau_{j'},s_{k',l'}).$$
By Assumption {\bf A2}, we have for any $\varepsilon>0$ and sufficiently large $T$,
$\vartheta(\varepsilon)<1.$
Note that $$\varrho(\tau_{j},s_{k,l};\tau_{j'},s_{k',l'})=r(\tau_{j},s_{k,l};\tau_{j'},s_{k',l'})+(1-r(\tau_{j},s_{k,l};\tau_{j'},s_{k',l'}))\rho(T)
\thicksim r(\tau_{j},s_{k,l};\tau_{j'},s_{k',l'})$$
uniformly for $(\tau_{j},s_{k,l}), (\tau_{j'},s_{k',l'})$ in the same $I_{i}$. Split the first sum into two parts as
\begin{eqnarray}\label{tt2}
\sum_{(\tau_{j},s_{k,l}),(\tau_{j'},s_{k',l'})\in I_{i}, 1\leq i\leq K_{T},\atop |s_{k,l}-s_{k',l'}|< \varepsilon, |\tau_{j}-\tau_{j'}|< \varepsilon}M(u)
+\sum_{(\tau_{j},s_{k,l}),(\tau_{j'},s_{k',l'})\in I_{i}, 1\leq i\leq K_{T},\atop \{|s_{k,l}-s_{k',l'}|\geq \varepsilon\}\cup\{|\tau_{j}-\tau_{j'}|\geq \varepsilon\} }M(u)=: S_{T,1}+S_{T,2}
\end{eqnarray}
Assumption {\bf A2} implies that for all $|s_{k,l}-s_{k',l'}|< \varepsilon<2^{-1/\alpha}$ and $|\tau_{j}-\tau_{j'}|< \varepsilon<2^{-1/\alpha}$
$$1-r(\tau_{j},s_{k,l};\tau_{j'},s_{k',l'})\leq C[|s_{k,l}-s_{k',l'}|^{\alpha}+|s_{k,l}-s_{k',l'}+\tau_{j}-\tau_{j'}|^{\alpha}].$$
Thus, we have
\begin{eqnarray*}
S_{T,1}&\leq&C\sum_{(\tau_{j},s_{k,l}),(\tau_{j'},s_{k',l'})\in I_{i}, 1\leq i\leq K_{T},\atop |s_{k,l}-s_{k',l'}|< \varepsilon, |\tau_{j}-\tau_{j'}|< \varepsilon}
\frac{|1-r(\tau_{j},s_{k,l};\tau_{j'},s_{k',l'})|}{\sqrt{1-\varrho(\tau_{j},s_{k,l};\tau_{j'},s_{k',l'})}}\rho(T)
\exp\left(-\frac{u^{2}}{1+r(\tau_{j},s_{k,l};\tau_{j'},s_{k',l'})}\right)\\
&\leq & C\sum_{(\tau_{j},s_{k,l}),(\tau_{j'},s_{k',l'})\in I_{i}, 1\leq i\leq K_{T},\atop |s_{k,l}-s_{k',l'}|< \varepsilon, |\tau_{j}-\tau_{j'}|< \varepsilon}
\sqrt{1-r(\tau_{j},s_{k,l};\tau_{j'},s_{k',l'})}\rho(T)
\exp\left(-\frac{u^{2}}{2}\right)\\
&\leq&C\sum_{(\tau_{j},s_{k,l}),(\tau_{j'},s_{k',l'})\in I_{i}, 1\leq i\leq K_{T},\atop |s_{k,l}-s_{k',l'}|< \varepsilon, |\tau_{j}-\tau_{j'}|< \varepsilon}
[|s_{k,l}-s_{k',l'}|^{\alpha}+|s_{k,l}-s_{k',l'}+\tau_{j}-\tau_{j'}|^{\alpha}]^{1/2}\rho(T)
\exp\left(-\frac{u^{2}}{2}\right).
\end{eqnarray*}
Note that
$$Tw(u)=O(1),\ \ u\rightarrow\infty,$$
which implies
$$u^{2}=2\ln T+(\frac{4}{\alpha}-1)\ln\ln T+O(1).$$
Since there are $C(1/q^{2})$ combinations of two points
$\tau_{j}, \tau_{j'}\in \cup_{i} I_{i}$ and $T\varepsilon/q^{2}$ combinations of two points
$s_{k,l}, s_{k',l'}\in \cup_{i} I_{i}$.
We have  (recall $\rho(T)=r/\ln T$ and $q=du^{-2/\alpha}$)
\begin{eqnarray*}
S_{T,1}&\leq& C Tq^{-2}\rho(T)\exp\left(-\frac{u^{2}}{2}\right)\sum_{0<kq<\varepsilon}\sum_{0<lq<\varepsilon}[|kq|^{\alpha}+|kq+lq|^{\alpha}]^{1/2}\\
&\leq& C Tq^{-2}\rho(T)\exp\left(-\frac{u^{2}}{2}\right)\\
&\leq&C(\ln T)^{-1/2}.
\end{eqnarray*}
Thus,  $S_{T,1}\rightarrow0$ as $T\rightarrow\infty$.
For the second sum, similarly, we have
\begin{eqnarray*}
S_{T,2}&\leq&C\sum_{(\tau_{j},s_{k,l}),(\tau_{j'},s_{k',l'})\in I_{i}, 1\leq i\leq K_{T},\atop \{|s_{k,l}-s_{k',l'}|\geq \varepsilon\}\cup\{|\tau_{j}-\tau_{j'}|\geq \varepsilon\} }
|r(\tau_{j},s_{k,l};\tau_{j'},s_{k',l'})-\varrho(\tau_{j},s_{k,l};\tau_{j'},s_{k',l'})|\\
&&\times
\exp\left(-\frac{u^{2}}{1+\varpi(\tau_{j},s_{k,l};\tau_{j'},s_{k',l'})}\right)\\
&\leq & C\sum_{(\tau_{j},s_{k,l}),(\tau_{j'},s_{k',l'})\in I_{i}, 1\leq i\leq K_{T},\atop \{|s_{k,l}-s_{k',l'}|\geq \varepsilon\}\cup\{|\tau_{j}-\tau_{j'}|\geq \varepsilon\} }
\exp\left(-\frac{u^{2}}{1+\vartheta(\varepsilon)}\right)\\
&\leq&C Tq^{-4}\exp\left(-\frac{u^{2}}{1+\vartheta(\varepsilon)}\right)\\
&\leq&C T^{-\frac{1-\vartheta(\varepsilon)}{1+\vartheta(\varepsilon)}}(\ln T)^{4/\alpha}.
\end{eqnarray*}
Since $\vartheta(\varepsilon)<1$, we get $S_{T,2}\rightarrow0$ as $T\rightarrow\infty$.

We continue to estimate the upper bound for the second sum in (\ref{tt1}).
\zT{Note that}, in this case, $|s_{k,l}-s_{k',l'}|\geq \delta$ and $\varrho(\tau_{j},s_{k,l},\tau_{j'},s_{k',l'})=\rho(T)=r/\ln T$,
\zT{then} by Assumptions {\bf A2} and {\bf A3},
$$\varpi(\tau_{j},s_{k,l};\tau_{j}',s_{k',l'})=\max\{|r(\tau_{j},s_{k,l};\tau_{j'},s_{k',l'})|,|\varrho(\tau_{j},s_{k,l};\tau_{j'},s_{k',l'})|\}<\rho<1$$
 for $|s_{k,l}-s_{k',l'}|\geq \delta$. Set
$\beta<(1-\rho)/(1+\rho)$ and
split the second sum into two parts as
\begin{eqnarray}\label{tt3}
\sum_{(\tau_{j},s_{k,l})\in I_{i},(\tau_{j'},s_{k',l'})\in I_{i'},\atop 1\leq i\neq i'\leq K_{T}, |s_{k,l}-s_{k',l'}|< T^{\beta}}M(u)
+\sum_{(\tau_{j},s_{k,l})\in I_{i},(\tau_{j'},s_{k',l'})\in I_{i'}, \atop 1\leq i\neq i'\leq K_{T}, |s_{k,l}-s_{k',l'}|\geq T^{\beta}}M(u)=: S_{T,3}+S_{T,4}.
\end{eqnarray}
For the first
sum, there are $T^{1+\beta}/q^{2}$ combinations of two points
$s_{k,l}, s_{k',l'}\in \cup_{k} I_{k}$. Together with the $\tau_{j}$
combinations, there are $(T^{1+\beta}/q^{2})(1/q^{2})$ terms in
the sum $S_{T,3}$.
Thus, for $S_{T,3}$ we have
\begin{eqnarray*}
S_{T,3}&\leq&C \sum_{(\tau_{j},s_{k,l})\in I_{i},(\tau_{j'},s_{k',l'})\in I_{i'},\atop 1\leq i\neq i'\leq K_{T}, |s_{k,l}-s_{k',l'}|< T^{\beta}}|r(\tau_{j},s_{k,l};\tau_{j'},s_{k',l'})-\varrho(\tau_{j},s_{k,l};\tau_{j'},s_{k',l'})|\\
&&\times
\exp\left(-\frac{u^{2}}{1+\varpi(\tau_{j},s_{k,l},\tau_{j'},s_{k',l'})}\right)\\
&\leq&C \frac{T^{1+\beta}}{q^{4}}\exp\left(-\frac{u^{2}}{1+\rho}\right)\\
&=&CT^{1+\beta-2/(1+\rho)}(\ln T)^{4/\alpha}\rightarrow0
\end{eqnarray*}
 $T \to \IF$\zT{, because
of $1+\beta<2/(1+\rho)$, which is due to} the choice of
$\beta$, \zT{and by} using $q=du^{-2/\alpha}$.\\
For the second sum $S_{T,4}$ with $|s_{k,l}-s_{k',l'}|\geq
T^{\beta}$, by Assumption {\bf A3 }we have
$$\sup_{|s_{k,l}-s_{k',l'}|\geq T^{\beta}}\varpi(\tau_{j},s_{k,l};\tau_{j}',s_{k',l'})\leq C(\ln(T^{\beta}))^{-1}.$$
So, for $|s_{k,l}-s_{k',l'}|\geq T^{\beta}$, we have
\begin{eqnarray*}
&&\frac{T^{2}}{q^{4}\ln T}\exp\left(-\frac{u^{2}}{1+\varpi(\tau_{j},s_{k,l};\tau_{j'},s_{k',l'})}\right)
\leq\frac{T^{2}}{q^{4}\ln T}\exp\left(-\frac{u^{2}}{1+C/\ln T^{\beta}}\right)=O(1)
\end{eqnarray*}
Hence, we have
\begin{eqnarray}
\label{sum1}
S_{T,4}&\leq& C\sum_{(\tau_{j},s_{k,l})\in I_{i},(\tau_{j'},s_{k',l'})\in I_{i'}, \atop 1\leq i\neq i'\leq K_{T}, |s_{k,l}-s_{k',l'}|\geq T^{\beta}}|r(\tau_{j},s_{k,l};\tau_{j'},s_{k',l'})-\varrho(\tau_{j},s_{k,l};\tau_{j'},s_{k',l'})|\nonumber\\
&&\times\exp\left(-\frac{u^{2}}{1+\varpi(\tau_{j},s_{k,l};\tau_{j'},s_{k',l'})}\right)\nonumber\\
&=& C\frac{q^{4}\ln T}{T^{2}}\sum_{(\tau_{j},s_{k,l})\in I_{i},(\tau_{j'},s_{k',l'})\in I_{i'}, \atop 1\leq i\neq i'\leq K_{T}, |s_{k,l}-s_{k',l'}|\geq T^{\beta}}|r(\tau_{j},s_{k,l};\tau_{j'},s_{k',l'})-\rho(T)|\nonumber\\
&&\times\frac{T^{2}}{q^{4}\ln T}\exp\left(-\frac{u^{2}}{1+\varpi(\tau_{j},s_{k,l};\tau_{j'},s_{k',l'})}\right)\nonumber\\
&\leq& C\frac{q^{4}\ln T}{T^{2}}\sum_{(\tau_{j},s_{k,l})\in I_{i},(\tau_{j'},s_{k',l'})\in I_{i'}, \atop 1\leq i\neq i'\leq K_{T}, |s_{k,l}-s_{k',l'}|\geq T^{\beta}}|r(\tau_{j},s_{k,l};\tau_{j'},s_{k',l'})-\rho(T)|\nonumber\\
&\leq& C\frac{q^{4}}{T^{2}}\sum_{(\tau_{j},s_{k,l})\in I_{i},(\tau_{j'},s_{k',l'})\in I_{i'}, \atop 1\leq i\neq i'\leq K_{T}, |s_{k,l}-s_{k',l'}|\geq T^{\beta}}|r(\tau_{j},s_{k,l};\tau_{j'},s_{k',l'})\ln(|s_{k,l}-s_{k',l'}|)-r|\nonumber\\
&&+Cr\frac{q^{4}}{T^{2}}\sum_{(\tau_{j},s_{k,l})\in I_{i},(\tau_{j'},s_{k',l'})\in I_{i'}, \atop 1\leq i\neq i'\leq K_{T}, |s_{k,l}-s_{k',l'}|\geq T^{\beta}}
\left|1-\frac{\ln T}{\ln(|s_{k,l}-s_{k',l'}|)}\right|
\end{eqnarray}
By Assumption {\bf A3}, the first term on the  right-hand side of (\ref{sum1}) tends to $0$ as $T\rightarrow\infty$.
Furthermore, the second term of the right-hand-side of (\ref{sum1}) also tends to 0 by an integral estimate as follows. We have
\begin{eqnarray*}
&&Cr\frac{q^{4}}{T^{2}}\sum_{(\tau_{j},s_{k,l})\in I_{i},(\tau_{j'},s_{k',l'})\in I_{i'}, \atop 1\leq i\neq i'\leq K_{T}, |s_{k,l}-s_{k',l'}|\geq T^{\beta}}
\left|1-\frac{\ln T}{\ln(|s_{k,l}-s_{k',l'}|)}\right|\\
&&\leq Cr\frac{1}{\beta\ln T}\frac{q^{2}}{T^{2}}\sum_{T^{\beta}<|kq-lq|\leq T}\left|\ln\left(\frac{|kq-lq|}{T}\right)\right|\\
&&\leq Cr\frac{1}{\beta\ln T}\int_{0}^{1}\int_{0}^{1}|\ln(x-y)|dxdy,
\end{eqnarray*}
which tends to $0$, as $T\rightarrow\infty$.
This completes the proof of the lemma.
\hfill$\Box$

\textbf{Proof of Theorem 2.2:}
\zT{First, by \nelem{Lemma 3.2}, \nelem{Lemma 3.3} and
\nelem{Lemma 3.4}, it holds that} as $T \to \IF$
\begin{eqnarray*}
P\left( \max_{(\tau,s)\in[a,b]\times[0,T]} \ZHTY \leq u\right)
&\sim&P\left(\max_{(\tau,s)\in\cup_{k} I_{k}} \ZHTY \leq u\right)\\
&\sim&P\left(\max_{(\tau,s)\in\cup_{k} I_{k}\cap \mathcal{R}} \ZHTY \leq u\right)\\
&\sim&P\left(\max_{(\tau,s)\in \cup_{k}I_{k}\cap \mathcal{R}} \YHTY \leq u\right)\\
&\sim&P\left(\max_{(\tau,s)\in \cup_{k}I_{k}} \YHTY \leq u\right)
\end{eqnarray*}
uniformly for $d>0$.
\zT{By the definition of $\YHTY $,} we have
\BQNY
P\left(\max_{(\tau,s)\in \cup_{k}I_{k}} \YHTY \leq u\right)
&=&P\left(\max_{(\tau,s)\in \cup_{k}I_{k}}(1-\rho(T))^{1/2}\xi(\tau,s)/\sigma(\tau)+(\rho(T))^{1/2}\mathcal{N}\leq u\right)\\
&=& \int_{-\infty}^{+\infty}P\left(\max_{(\tau,s)\in \cup_{k}I_{k}}\xi(\tau,s)/\sigma(\tau)\leq \frac{u-(\rho(T))^{1/2}z}{(1-\rho(T))^{1/2}}\right)d\Phi(z)\\
&=&\int_{-\infty}^{+\infty}\prod_{i=1}^{K_{T}}P\left(\max_{(\tau,s)\in I_{k}}\ZHTY\leq \frac{u-(\rho(T))^{1/2}z}{(1-\rho(T))^{1/2}}\right)d\Phi(z)\\
&=&\int_{-\infty}^{+\infty}K_{T}P\left(\max_{(\tau,s)\in I_{1}}\ZHTY\leq \frac{u-(\rho(T))^{1/2}z}{(1-\rho(T))^{1/2}}\right)d\Phi(z),
\EQNY
where $\Phi(\cdot)$ denotes the distribution function of a standard Gaussian random variable.
Recall that $u=u_{T}(x)=a_{T}^{-1}x+b_{T}$.
By some standard computations, we have
$$u_{z}:=\frac{u-(\rho(T))^{1/2}z}{(1-\rho(T))^{1/2}}=u+\frac{r+\sqrt{2r}z}{u}+o(u^{-1}),$$
as $T\rightarrow\infty$. Thus, by Theorem 2.1 and the definition on $K_{T}$, we have
\BQNY
P\left( \max_{(\tau,s)\in[a,b]\times[0,T]} \ZHTY \leq u\right)
&\sim&\int_{-\infty}^{+\infty}\exp\left(-K_{T}P\left(\max_{(\tau,s)\in I_{1}} \ZHTY > u_{z}\right)\right)d\Phi(z)\\
&\sim&\int_{-\infty}^{+\infty}\exp\left(-K_{T}(L-\delta)w(u_{z})\right)d\Phi(z)\\
&\rightarrow&\int_{-\infty}^{+\infty} \exp(-e^{-x-r+\sqrt{2r}z}) d\Phi(z)
\EQNY
as $\delta\downarrow0,\ \ T\to \IF$.
Therefore, the claim follows. \hfill$\Box$

\subsection{Proof of Propositions 3.1 and 3.2}

\textbf{Proof of Proposition 3.1:} In order to prove this proposition, it suffices to check that Assumptions {\bf A1-A3} hold.
For the stationary Gaussian process $X(t)$ with correlation function $r_{X}(t)$, it easy to see that
$$E(X(s+\tau)-X(\tau))^{2}=2(1-r_{X}(\tau)).$$
Thus, Assumption {\bf A1} holds with $\sigma(\tau)=(1-r_{X}(\tau))^{1/2}$.\\
It holds for the correlation function $r_{Y}(\tau,s,\tau',s')$ of $Y(\tau,s)$ that
$$r_{Y}(\tau,s;\tau+\vartriangle_{\tau},s+\vartriangle_{s})=\frac{r_{X}(|\vartriangle_{\tau}+\vartriangle_{s}|)
-r_{X}(|\tau-\vartriangle_{s}|)-r_{X}(|\vartriangle_{\tau}+\tau+\vartriangle_{s}|)+r_{X}(|\vartriangle_{s}|)}
{\sqrt{2(1-r_{X}(\tau))2(1-r_{X}(\tau+\vartriangle_{\tau}))}}.$$
By Assumption {\bf B1}, we can see that
$$r_{Y}(\tau,s;\tau+\vartriangle_{\tau},s+\vartriangle_{s})=1-\frac{a_{1}|\vartriangle_{\tau}+\vartriangle_{s}|^{\alpha}+a_{1}|\vartriangle_{s}|^{\alpha}}{2(1-r_{X}(\tau))}
(1+o(1))$$
as $\vartriangle_{\tau}\rightarrow0$ and $\vartriangle_{s}\rightarrow0$. Thus, Assumption {\bf A2} holds with $g(\tau)=a_{1}/(2(1-r_{X}(\tau)))$.\\
Since $r_{X}(t)$ is twice continuously differentiable in $(0,\infty)$, we have
$$|r_{X}(|\tau+s-\tau'-s'|)-r_{X}(|s+s'+\tau'|)-r_{X}(|s-s'-\tau'|)+r_{X}(|s-s'|)|\leq C \ddot{r}_{X}(|s-s'|)$$
for $\tau,\tau'\in[a,b]$ as $|s-s'|\rightarrow\infty$,
which together with Assumption {\bf B2} implies
Assumption {\bf A3}. $\Box$

\textbf{Proof of Proposition 3.2:} As for the proof of Proposition 3.1, we check that Assumptions ${\bf A1-A3}$ hold. Using the stationarity of the increments of $X(t)$ and Assumption {\bf C1}, it follows that Assumption {\bf A1} holds with $\sigma(\tau)=\sigma_{X}(\tau)=(E(X(\tau)))^{1/2}$.\\
From the stationarity of the increments of $X(t)$ we have
$$Cov(X(t),X(s))=\frac{1}{2}[\sigma^{2}_{X}(t)+\sigma^{2}_{X}(s)-\sigma^{2}_{X}(|t-s|)],$$
which implies the correlation function of $Y(\tau,s)$ equals
$$r_{Y}(\tau,s;\tau',s')=\frac{1}{2\sigma_{X}(\tau)\sigma_{X}(\tau')}[\sigma_{X}^{2}(|\tau+s-\tau'-s'|)-\sigma_{X}^{2}(|s-s'+\tau|)
+\sigma_{X}^{2}(|s-s'+\tau'|)-\sigma_{X}^{2}(|s-s'|)].$$
Thus, it follows from Assumption {\bf C1} that
$$r_{Z}(\tau,s;\tau+\vartriangle_{\tau},s+\vartriangle_{s})=1-\frac{a_{2}}{2\sigma_{X}^{2}(\tau)}[|\vartriangle_{\tau}+\vartriangle_{s}|^{\alpha}+|\vartriangle_{s}|^{\alpha}](1+o(1)),$$
as $\vartriangle_{\tau}\rightarrow 0$ and $\vartriangle_{s}\rightarrow 0$ and Assumption {\bf A2} holds with $g(\tau)=a_{2}/(2\sigma_{X}^{2}(\tau))$. \\
By Taylor expansions, it is straightforward to check that
$$|r_{Y}(s,t;s',t')|\leq C\ddot{\sigma}_{X}^{2}(|s-s'|)$$
for $\tau,\tau'\in[a,b]$ as $|s-s'|\rightarrow\infty$, which combined with Assumption {\bf C2} implies Assumption {\bf A3}.
$\Box$

\section{Appendix}

In this subsection, we \zT{make a brief review of} some
important results and definitions for locally stationary Gaussian
fields. These results are obtained by Piterbarg and Prisyazhnyuk (1978), Mikhaleva and Piterbarg
(1996), Chan and Lai (2006)  and Kabluchko (2011). These results can also be found in the Appendix of Tan and Yang (2015).

 \BD \label{Definition 4.1} A function $f: \mathbb{R}^{d}\rightarrow \mathbb{R}$ is called homogeneous of order
 $\alpha>0$,
 if for each $\mathbf{s}\in \mathbb{R}^{d}$ and $\lambda\in
 \mathbb{R}$,
$$f(\lambda\mathbf{s})=|\lambda|^{\alpha}f(\mathbf{s}).$$
\ED

Let $\mathcal{E}(\alpha)$ be the set of all continuous homogeneous
functions of order $\alpha$. For $f\in \mathcal{E}(\alpha)$,
\zT{define} $||f||= \max_{||\mathbf{t}||_{2}=1}f(\mathbf{t})$.
Denote by $\mathcal{E}^{+}(\alpha)$ the cone of all strictly
positive functions in $\mathcal{E}(\alpha)$.

 \BD \label{Definition 4.2} Let $\{\xi(\mathbf{t}),\mathbf{t}\in D\}$ be a centered Gaussian field with constant variance $1$
 defined on some domain $D\subset \mathbb{R}^{d}$.
 Let $r(\mathbf{t}_{1},\mathbf{t}_{2})=\mathbb{E}\{\xi(\mathbf{t}_{1})\xi(\mathbf{t}_{2})\}$ be the covariance function
of $\xi$ and suppose that it satisfies
$r(\mathbf{t}_{1},\mathbf{t}_{2})=1\Leftrightarrow
\mathbf{t}_{1}=\mathbf{t}_{2}$. The Gaussian field $\xi$ is called
locally stationary with index $\alpha\in(0,2]$, if for each
$\mathbf{t}\in D$, there exits a continuous function
$C_{\mathbf{t}}\in \mathcal{E}^{+}(\alpha)$, such that
$$
\frac{1-r(\mathbf{t},\mathbf{t}+\mathbf{s})}{C_{\mathbf{t}}(\mathbf{s})} \to 1, \quad ||\mathbf{s}||_{2}\rightarrow0$$
holds uniformly on compacts, and the map $C_{\cdot}: D\rightarrow \mathcal{E}^{+}(\alpha)$ from $\mathbf{t}$ to $C_{\mathbf{t}}$ is continuous.
\ED

The collection of homogeneous $C_{\mathbf{t}}$ is referred to as the local structure of the field $\xi$.

We also need the following definition.

 \BD \label{Definition 4.3}Let $\{\xi(\mathbf{t}),\mathbf{t}\in D\}$ be a centered Gaussian field defined on some domain $D\subset \mathbb{R}^{d}$.
 Suppose that $\xi$ is locally stationary with index $\alpha$ and local structure
$C_{\mathbf{t}}(\mathbf{s})$. For each $\mathbf{t}\in D$, let $\{Y_{\mathbf{t}}(\mathbf{s}), \mathbf{s}\in \mathbb{R}^{d}\}$ be a
Gaussian field such that
$$\mathbb{E}\{Y_{\mathbf{t}}(\mathbf{s})\}=-C_{\mathbf{t}}(\mathbf{s}), \quad Cov(Y_{\mathbf{t}}(\mathbf{s}_{1}),Y_{\mathbf{t}}(\mathbf{s}_{2}))
=C_{\mathbf{t}}(\mathbf{s}_{1})+C_{\mathbf{t}}(\mathbf{s}_{2})-C_{\mathbf{t}}(\mathbf{s}_{1}-\mathbf{s}_{2}), \quad \mathbf{s},\mathbf{s}_1,
\mathbf{s}_2\inr^d.
$$
Then, $Y_{\mathbf{t}}$ is called the tangent field of $\xi$ at the
point $\mathbf{t}$ conditional on $\xi(\mathbf{t})=\infty$, and
$$h(\mathbf{t})=\lim_{\mathcal{T}\rightarrow\infty}\frac{1}{\mathcal{T}^{d}}\mathbb{E}
\left\{\exp\left(\max_{\mathbf{s}\in[0,\mathcal{T}]^{d}} Y_{\mathbf{t}}(\mathbf{s})\right)\right\}$$
is called the high excursion intensity of the field $\xi$.
\ED

Chan and Lai (2006) \zT{showed} that $h(\mathbf{t})\in (0,\infty)$
exists and is continuous in $\mathbf{t}$; the \zT{following} theorem
therein determines the asymptotic behavior of the high excursion
probability of a locally stationary Gaussian field.

\BT \label{Theorem 4.1} Let $\{\xi(\mathbf{t}),\mathbf{t}\in D\}$ be
a centered Gaussian field defined on some domain $D\subset
\mathbb{R}^{d}$. Suppose that $\xi$ is locally stationary with index
$\alpha$ and local structure $C_{\mathbf{t}}(\mathbf{s})$. Then, for
any compact, set $K\subset D$ with positive Jordan measure
$$P\left( \max_{\mathbf{t}\in K}\xi(\mathbf{t})>u\right)= \left(\int_{K} h(\mathbf{t})d\mathbf{t}\right)u^{\frac{2d}{\alpha}}\Psi(u)
(1+o(1)), \quad u \to \IF,$$
with  $h(\mathbf{t}): D\rightarrow(0,\infty)$ the high excursion intensity of $\xi$.
\ET

The following result, which describes the asymptotic behavior of the high excursion probability over a finite grid with mesh size going to $0$,
plays an important role in our proof. It's proof can be found in Kabluchko (2011).
\BT \label{Theorem 4.2} Suppose that the conditions of Theorem \ref{Theorem 4.1} are satisfied.
Let $q=q(u)=du^{-2/\alpha}$ Then, for
any compact set $K\subset D$ with positive Jordan measure
$$P\left( \max_{\mathbf{t}\in K\cap q\mathbb{Z}^{d}}\xi(\mathbf{t})>u\right)= \left(\int_{K} h_{d}(\mathbf{t})d\mathbf{t}\right)u^{\frac{2d}{\alpha}}\Psi(u)
(1+o(1)), \quad u \to \IF,$$
where
$$h_{d}(\mathbf{t})=\lim_{\mathcal{T}\rightarrow\infty}\frac{1}{\mathcal{T}^{d}}\mathbb{E}
\left\{\exp\left(\max_{\mathbf{s}\in[0,\mathcal{T}]^{d}\cap d\mathbb{Z}^{d}} Y_{\mathbf{t}}(\mathbf{s})\right)\right\}.$$
Furthermore, $\lim_{d\downarrow0}h_{d}(\mathbf{t})=h(\mathbf{t})$, where $h(\mathbf{t})$ is the high excursion intensity of $X$.
\ET

\textbf{Acknowledgments.} \zT{The authors would like to thank the
referee and the editor for the thorough reading and valuable
suggestions, which greatly improve the original results of the paper.}

\end{document}